\documentclass[twocolumn]{autart}    

\usepackage[T1]{fontenc} 
\usepackage{lmodern}  

\usepackage{graphicx}          
\usepackage{cite}
\usepackage{amsmath,amssymb}
\usepackage{algorithmic}
\usepackage{graphicx, subfigure}
\usepackage{textcomp}
\usepackage{xcolor}
\usepackage{enumitem}
\usepackage{xfrac} 
\usepackage{mleftright}
\usepackage{IEEEtrantools}

\newcommand{\zz}{z}

\newcommand{\pp}{\psi}

\newcommand{\xs}{x^\star}

\newcommand{\zs}{\zz^\star}

\newcommand{\pps}{\pp^\star}

\newcommand{\xd}{\dot{x}}
\newcommand{\zd}{\dot{\zz}}

\newcommand{\xdd}{\ddot{x}}
\newcommand{\zdd}{\ddot{\zz}}

\newcommand{\xt}{\tilde{x}}

\newcommand{\Lt}{\tilde{L}}
\newcommand{\mt}{\widetilde{m}}
\newcommand{\pt}{\tilde{p}}

\newcommand{\xp}{x^\prime}

\renewcommand{\tt}{\theta}

\newcommand{\cF}{\enma{\mathcal F}}
\newcommand{\cE}{\enma{\mathcal E}}

\newcommand{\cM}{\enma{\mathcal M}}
\newcommand{\cC}{\enma{\mathcal C}}

\newcommand{\cD}{\enma{\mathcal D}}

\newcommand{\cFm}{\cF_\mu}
\newcommand{\cDm}{\cD_\mu}

\makeatletter
\newcommand*\rel@kern[1]{\kern#1\dimexpr\macc@kerna}
\newcommand*\widebar[1]{%
	\begingroup
	\def\mathaccent##1##2{%
		\rel@kern{0.8}%
		\overline{\rel@kern{-0.8}\macc@nucleus\rel@kern{0.2}}%
		\rel@kern{-0.2}%
	}%
	\macc@depth\@ne
	\let\math@bgroup\@empty \let\math@egroup\macc@set@skewchar
	\mathsurround\z@ \frozen@everymath{\mathgroup\macc@group\relax}%
	\macc@set@skewchar\relax
	\let\mathaccentV\macc@nested@a
	\macc@nested@a\relax111{#1}%
	\endgroup
}

\DeclareMathOperator*{\argmin}{argmin}
\DeclareMathOperator*{\argmax}{argmax}

\DeclareMathOperator*{\minimize}{minimize}
\DeclareMathOperator*{\maximize}{maximize}

\newcommand{\reals}[1]{\mathbb{R}^{#1}}

\DeclareMathOperator{\prox}{\mathbf{prox}}

\newcommand{\norm}[1]{\| #1 \|}    
\newcommand{\inner}[2]{\langle #1, #2\rangle}    

\newcommand{\set}[1]{\{ #1 \}}
\newcommand{\enma}[1]{\ensuremath{#1}}

\newcommand{\non}{\nonumber}

\newcommand{\beq}{\begin{equation}}
\newcommand{\eeq}{\end{equation}}
\newcommand{\ba}{\begin{array}}
\newcommand{\ea}{\end{array}}
\newcommand{\bseq}{\begin{subequations}}
\newcommand{\eseq}{\end{subequations}}

\newcommand{\DefinedAs}[0]{\mathrel{\mathop:}=}

\newcommand{\rmd}{\mathrm{d}}
\newcommand{\rme}{\mathrm{e}}

\newcommand{\matbegin}{\left[}
\newcommand{\matend}{\right]}

\newcommand{\tbo}[2]{
		\matbegin \begin{array}{c}
				#1 \\ #2
			\end{array} \matend }
\newcommand{\thbo}[3]{
		\matbegin \begin{array}{c}
				#1 \\ #2 \\ #3
			\end{array} \matend }

\newcommand{\tbt}[4]{
		\matbegin \begin{array}{cc}
				#1 & #2 \\ #3 & #4
			\end{array} \matend }

\newcommand{\thbth}[9]{
		\matbegin \begin{array}{ccc}
				#1 & #2 & #3 \\
				#4 & #5 & #6 \\
				#7 & #8 & #9
			\end{array}\matend}

\newtheorem{theorem}{Theorem}

\newtheorem{lemma}{Lemma}
\newtheorem{remark}{Remark}
\newtheorem{assumption}{Assumption}

\newcommand{\nvs}{0.1cm}

\newcommand{\asp}[1]{ #1 }

\begin{document}

\begin{frontmatter}
\runtitle{Accelerated forward-backward and Douglas-Rachford splitting dynamics}  

\title{Accelerated forward-backward and \\[0.15cm]Douglas-Rachford splitting dynamics}
	\vspace*{-4ex}
\author{Ibrahim K.\ Ozaslan and~}\ead{ozaslan@usc.edu}    
\author{Mihailo R.\ Jovanovi\'c}\ead{mihailo@usc.edu}              
\address{Ming Hsieh Department of Electrical and Computer Engineering, University of Southern California, Los Angeles, CA 90089}  


\begin{abstract}                          
	We examine convergence properties of continuous-time variants of accelerated Forward-Backward (FB) and Douglas-Rachford (DR) splitting algorithms for nonsmooth composite optimization problems. When the objective function is given by the sum of a quadratic and a nonsmooth term, we establish accelerated sublinear and exponential convergence rates for convex and strongly convex problems, respectively. Moreover, for FB splitting dynamics, we demonstrate that accelerated exponential convergence rate carries over to general strongly convex problems. In our Lyapunov-based analysis we exploit the variable-metric gradient interpretations of FB and DR splittings to obtain smooth Lyapunov functions that allow us to establish accelerated convergence rates. We provide computational experiments to demonstrate the merits and the \mbox{effectiveness of our analysis.} \\[-0.25cm]
\end{abstract}
\vspace*{-2ex}   

\end{frontmatter}

\section{Introduction}\label{sec.intro}
\vspace*{-1ex}

The convergence properties of gradient descent (GD),
\beq\label{eq.gd}
x^{k+1}  \,=\,  x^k  \,-\,  \widebar{\alpha}\nabla f(x^k)
\eeq
can be improved by adding a momentum term,
\beq\label{eq.acc_gd}
x^{k + 1} \,=\, x^k  \,+\, \widebar{\gamma}(x^{k} \,-\, x^{k - 1}) \,+\, \widebar{\alpha}\nabla f(x^k \, + \, \widebar{\beta} (x^k  \, - \,  x^{k-1})).
\eeq 
where $k$ is the iteration index, $\widebar{\alpha}$ is the stepsize, $\widebar{\beta}$ is the extrapolation parameter, and $\widebar{\gamma}$ is the damping coefficient. In particular, the Polyak's Heavy-Ball Method (HBM) is obtained by setting $\widebar{\beta} = 0$ in~\eqref{eq.acc_gd}~\cite{pol64} and Nesterov's Accelerated Algorithm (NAA) results from $\widebar{\gamma} = \widebar{\beta}$ in~\eqref{eq.acc_gd}~\cite{nes83}. For a convex objective function $f$: $\reals{n}\to\reals{}$ with an $L$-Lipschitz continuous gradient $\nabla f$, NAA with $(\widebar{\alpha}, \widebar{\beta}) = (1/L, (k-2)/(k+3))$ improves the sublinear convergence rate $O(1/ k)$ of GD to $O(1/k^2)$. Furthermore, for an $m$-strongly convex $f$, the linear convergence rate $O(\rme^{-k/\kappa})$ of GD~\cite[Section 2.1.5]{nes03} is improved to $O(\rme^{-k/\sqrt{\kappa}})$ with $(\widebar{\alpha}, \widebar{\beta}) = (1/L, (\sqrt{\kappa}-1)/(\sqrt{\kappa}+1) )$. In both cases, these rates are optimal in the sense that there exists a smooth convex function $f$ for which the rates cannot be improved~\cite[Theorems 2.1.7 and 2.1.13]{nes03}.   

A wide-spread use of accelerated algorithms has inspired efforts to understand the underlying mechanisms. In particular, a growing body of literature views optimization algorithms as continuous-time dynamical systems that allow familiar methods to be recovered upon proper discretization~\cite{frarobvid21, muejor19, dhikhojovTAC19, allcor23}. For example, the continuous-time equivalent of GD is given by the gradient flow dynamics, 
\beq\label{eq.gd_flow}
\xd  \,=\, -\alpha\nabla f(x)
\eeq
and celebrated Arrow-Hurwicz-Uzawa primal-dual gradient flow dynamics can be used to effectively solve a class of constrained optimization problems~\cite{ozajovCDC23}. A dynamical systems perspective of NAA offered in~\cite{suboycan14} was further exploited by leveraging a variational formulation based on Lagrangian and Hamiltonian frameworks from classical mechanics~\cite{wibwiljor16, wilrecjor21} and contraction theory~\cite{cisbul22}. Furthermore, NAA was obtained using the semi-implicit Euler discretization of the differential equation~\cite{muejor19},
\beq\label{eq.inertial}
\xdd \,+\, \gamma \xd \,+\, \alpha\nabla f(x \,+\, \beta\xd) \,=\, 0
\eeq   
where $\alpha$, $\beta$, and $\gamma$ are possibly time-varying algorithmic parameters that differ for convex and strongly convex problems. The following Lyapunov function, 
\beq\label{eq.lyap0}
V(x, \xd) \,=\, af(x)  \,+\, (1/2)\norm{bx \,+\, \xd}_2^2
\eeq
with positive parameters $a$ and $b$ was also utilized in~\cite{muejor19} to show that~\eqref{eq.inertial} achieves the optimal convergence rate of discrete-time algorithms for both convex and strongly convex problems with proper choices of ($\alpha,\beta,\gamma$). Moreover, for $\beta + \gamma = 1$, the discretization of~\eqref{eq.inertial} preserves the accelerated convergence rate. This provides a system-theoretic characterization of acceleration and helps bridging the gap between discrete and continuous-time analyses.

\textbf{Our contributions.} We generalize accelerated dynamics~\eqref{eq.inertial} to nonsmooth composite optimization problems in which the objective function is given by the sum of a smooth and a nonsmooth term. Our approach is motivated by~\cite{mogjovAUT21} which obtains the Forward-Backward (FB) and Douglas-Rachford (DR) splitting dynamics by replacing $\nabla f$ in~\eqref{eq.gd_flow} with the generalized gradient map associated with the nonsmooth objective function. Based on the relation between gradient flow dynamics~\eqref{eq.gd_flow} and its accelerated variant~\eqref{eq.inertial}, we introduce accelerated FB and DR splitting dynamics by replacing gradient $\nabla f$ in~\eqref{eq.inertial} with the generalized gradient map. We also utilize a Lyapunov-based approach to establish accelerated sublinear and linear convergence rates that match the optimal discrete-time rates for convex and strongly convex problems, respectively. Our Lyapunov function is similar to~\eqref{eq.lyap0} with two important differences: (i) the Euclidean distance is replaced with a suitable weighted norm; and (ii) the FB or DR envelopes are used instead of the objective function. The latter choice is motivated by the variable-metric gradient method interpretation of associated splittings~\cite{patstebem14, patstebem14a}. We also note that, in our analysis, the condition $\beta + \gamma = 1$ is satisfied, which is an integral requirement for discretization analysis of~\cite{muejor19} to be applicable to the proposed accelerated dynamics.

Although accelerated rates for the FB and DR splittings are known in discrete-time, the continuous-time analysis in the problem settings that we examine has not been done before. An accelerated sublinear convergence rate matching our result has been proposed in~\cite{suboycan14}, but the approach relies on a dynamical system that is based on directional subgradients whose evaluation requires solving linear program
$
\argmax_{z \in \partial g(x)} z^T\dot{x}
$
at every time instant, where $g$ denotes the nonsmooth term in the objective function. Moreover, the existence and uniqueness of a solution to this dynamical system has been assumed without formally establishing it~\cite[Theorem 24]{suboycan14}. An accelerated sublinear convergence rate has also been proposed in~\cite{attchb15,attchbpeyred18} for a differential inclusion that models accelerated FB splitting, but no analysis was provided.

It is also worth mentioning that the FB and DR splittings are special cases of the primal-dual gradient flow dynamics associated with the proximal augmented Lagrangian~\cite{dhikhojovTAC19, mogjovAUT21, ozajovCDC22}. Thus, insights about acceleration of these splitting methods can also pave the way to accelerated primal-dual algorithms, which is of critical importance to the design of fast gradient-based methods for constrained optimization problems. We note that~\eqref{eq.inertial} was used as a motivation in~\cite{attchbfadria22, zenleiche23, hehufan22a} to accelerate primal-dual gradient flow dynamics. However, in contrast to~\cite{muejor19}, the algorithmic parameters are not bounded in time in these references. Thus, it is unclear if the primal-dual dynamics proposed in~\cite{attchbfadria22, zenleiche23, hehufan22a} allow rate-preserving discretization that leads to efficient first-order iterative schemes.

\textbf{Paper structure.} In Section~\ref{sec.back}, we introduce the problem setup and provide background material. In Section~\ref{sec.main}, we summarize our main findings, and in Section~\ref{sec.quad} we provide the proofs. In Section~\ref{sec.simulation}, we utilize computational experiments to demonstrate the merits of our analyses, and in Section~\ref{sec.conclusion}, we conclude the paper with remarks.

	\vspace*{-2ex}
\section{Problem formulation and preliminaries}\label{sec.back}
	\vspace*{-1ex}
	
In this section, we formulate the problem and provide background on the proximal operators; Moreau, FB, and DR envelopes; as well as FB and DR splittings. 

We consider convex optimization problem
\beq\label{eq.composite}
\minimize\limits_{x}~F(x) \;\DefinedAs\; f(x) \,+\, g(x)
\eeq
where $f$: $\reals{n}\to\reals{}$ is, unless noted otherwise, a continuously differentiable convex function with an $L$-Lipschitz continuous gradient $\nabla f$ and $g$: $\reals{n}\to\reals{}$ is a possibly non-differentiable closed proper convex function. 

	\vspace*{-1ex}
\subsection{Proximal operator and Moreau envelope}\label{sec.prox}
	\vspace*{-1ex}

The proximal operator of a proper lower semi-continuous convex function $g$ for a positive penalty parameter $\mu$ is the mapping defined as $\prox_{\mu g}(v) = (I + \mu\partial g)^{-1}(v),$ where $\partial g$ is a subdifferential of $g$ and $I$ is the identity map~\cite{parboy14}. Alternatively, $\prox_{\mu g}$ can be obtained as 
\beq\non
\prox_{\mu g}(v)  \,=\, \argmin\limits_{z}
\left( g(z) \,+\, \tfrac{1}{2\mu}\norm{z \,-\, v}_2^2 \right).
\eeq
The value function of this optimization problem determines the Moreau envelope associated with $g$,
\beq\label{eq.moreau}
\ba{rcl}
\!\!\!\cM_{\mu g}(v)  
&\asp{=}& g(\prox_{\mu g}(v)) \,+\, \frac{1}{2\mu}\norm{\prox_{\mu g}(v) \, - \, v}_2^2
\ea 
\eeq
which is, even for a non-differentiable $g$, a continuously differentiable function~\cite{parboy14}, 
\beq\label{eq.moreau_grad}
\nabla\cM_{\mu g}(v) \,=\, \tfrac{1}{\mu}(v \, - \, \prox_{\mu g}(v)) \, \in \, \partial g(\prox_{\mu g}(v)). 
\eeq

	\vspace*{-2ex}
\subsection{The FB and DR envelopes}\label{sec.fbe}
	\vspace*{-1ex}
	
Using Lipschitz continuity of $\nabla f$, we have the following upper bound on the objective function $F$ in~\eqref{eq.composite}~\cite{patstebem14},
\beq\non
\ba{l}
\!\! F(y) \leq  \; J(x,y) \; \DefinedAs \;
 \\ 
\!\! f(x) \,+\, g(y) \,+\, \tfrac{1}{2\mu}\norm{y \,-\, (x  \,-\, \mu\nabla f(x))}_2^2 \,-\, \frac{\mu}{2}\norm{\nabla f(x)}_2^2
\ea
\eeq
which holds for any $x, y\in\reals{n}$ and $\mu\in(0,1/L)$. Minimizing  $J(x,y)$ over $y$ yields the FB envelope of $F$~\cite{patstebem14}
\beq\non\label{eq.fb_env}
\ba{rrl}
\cF_{\mu}(x) & \DefinedAs& \minimize\limits_{y}~J(x,y)  
\\ 
& = &
f(x) \,+\, \cM_{\mu g}(x  \,-\,  \mu\nabla f(x)) \,-\, \tfrac{\mu}{2}\norm{\nabla f(x)}_2^2.
\ea 
\eeq
For a particular choice of a distance function (instead of the Euclidean norm), the FB envelope can also be interpreted as a generalized Moreau envelope of $F$~\cite{liupon17}. The DR envelope associated with~\eqref{eq.composite} is another useful value function resulting from evaluation of the FB envelope $\cF_{\mu}(x)$ at $x = \prox_{\mu f}(z)$ with $z\in\reals{n}$~\cite{patstebem14a, gisfae18}, i.e., 
\beq\non
\cDm(z) \,\DefinedAs\, \cFm(\prox_{\mu f}(z)).
\eeq
When $f$ is a twice continuously differentiable function with a bounded Hessian $\nabla^2f(x)$, the FB envelope is also continuously differentiable with gradient
\beq\label{eq.fb_grad}
\nabla \cF_{\mu}(x) \,=\, (I \,-\, \mu\nabla^2f(x))G_\mu(x)
\eeq
where $G_\mu$ is the generalized gradient map defined as
\beq\label{eq.gen_grad}
G_\mu(x)  \,=\, \tfrac{1}{\mu}(x \,-\, \prox_{\mu g}(x \,-\, \mu\nabla f(x))).
\eeq
Consequently, the DR envelope is also continuously differentiable and its gradient is determined by
\beq\label{eq.dr_grad}
\nabla\cDm(z)  \,=\, \nabla\prox_{\mu f}(z)\nabla\cFm(\prox_{\mu f}(z)).
\eeq
Here, $\nabla\prox_{\mu f}(z)$ can be obtained by differentiating~\eqref{eq.moreau_grad}
\beq\label{eq.prox_grad}
\nabla\prox_{\mu f}(z) \;\DefinedAs\; (I \,+\, \mu\nabla^2f(\prox_{\mu f}(z)))^{-1}
\eeq
Using~\eqref{eq.fb_grad} and~\eqref{eq.dr_grad}, $\nabla\cDm$ can also be expressed in terms of the generalized gradient map as
\beq\label{eq.dr_grad2}
\nabla\cDm(z)  \,=\,  (2\nabla\prox_{\mu f}(z)  \,-\, I)G_\mu(\prox_{\mu f}(z)).
\eeq
We note that for $\mu\in(0,~1/L)$, 
\beq\label{eq.env_equivalence}
\ba{c}
\min~F  \,=\,  \min ~\cFm  \,=\,  \min ~\cDm
\\
\argmin \;F  \,=\,  \argmin \; \cFm  \,=\,  \prox_{\mu f}(\argmin \; \cDm)
\ea
\eeq
which together with the smoothness of the envelope functions motivates studying acceleration of dynamics~\eqref{eq.fb_dyn}~and~\eqref{eq.dr_dyn} through the FB and DR envelopes~\cite{patstebem14}. 

	\vspace*{-1ex}
\subsection{The FB and DR splittings}\label{sec.dre}
		\vspace*{-1ex}
		
Proximal gradient algorithm~\cite{parboy14},
\beq\label{eq.fb_split}
x^{k+1} \,=\, x^k  \,-\, \widebar{\alpha} G_\mu(x^k)
\eeq
replaces $\nabla f$ in~\eqref{eq.gd} with the generalized gradient map~$G_\mu$, thereby extending gradient descent to composite optimization problem~\eqref{eq.composite}. Iterations~\eqref{eq.fb_split} are also known as the FB splitting and their continuous-time equivalent is determined by the proximal gradient flow dynamics~\cite{mogjovAUT21},
\beq\label{eq.fb_dyn}
\xd \,=\, -\alpha G_\mu (x).
\eeq
Furthermore, evaluation of the generalized gradient map at $x=\prox_{\mu f}(z)$ facilitates extension of FB splitting~\eqref{eq.fb_split} to composite optimization problem~\eqref{eq.composite} in which both $f$ and $g$ are allowed to be non-differentiable, 
\beq
\label{eq.dr_split0}
z^{k+1} \;\in\; z^k  \,-\, \widebar{\alpha} G_\mu(\prox_{\mu f}(z^k)).
\eeq
Using~\eqref{eq.gen_grad} and setting $\widebar{\alpha} = \mu$ in~\eqref{eq.dr_split0} yields,
\beq
	\label{eq.DRsplitting}
\ba{rcl}
z^{k+1}  
	& \asp{\in} &
	z^k
	\, - \, 
	\prox_{\mu f}(z) \; +
	\\
& & \prox_{\mu g}(\prox_{\mu f}(z) \,-\, \mu\partial f(\prox_{\mu f}(z))).
\ea
\eeq
Moreover,~\eqref{eq.moreau_grad} can be used to rewrite~\eqref{eq.DRsplitting} as,
\beq\label{eq.dr_split}
z^{k+1} \,=\, z^k  \,-\, \prox_{\mu f}(z) \; + \; \prox_{\mu g}(2\prox_{\mu f}(z) - z)
\eeq  
which is known as DR splitting~\cite{ryuboy16}. While~\eqref{eq.composite} can be solved via~\eqref{eq.dr_split} even for a non-differentiable $f$, when $f$ is continuously differentiable $G_\mu (\prox_{\mu f}(z))$ becomes a single-valued function and~\eqref{eq.dr_split0} can be used to obtain the DR splitting dynamics~\cite{mogjovAUT21},
	\beq
	\label{eq.dr_dyn}
	\zd \,=\, -\alpha G_\mu (\prox_{\mu f}(z)).
	\eeq
Lemma~\ref{lemma.sepideh_ges_fb} characterizes stability properties of proximal gradient flow~\eqref{eq.fb_dyn} and DR splitting dynamics~\eqref{eq.dr_dyn} under the following assumption.

\begin{assumption}\label{ass.1}
In~\eqref{eq.composite}, $f$ is an $m$-strongly convex function with an $L$-Lipschitz continuous gradient $\nabla f$ and $g$ is a proper closed convex function.
\end{assumption}

\begin{lemma}[Thms.~3 and 7,~\cite{mogjovAUT21}]\label{lemma.sepideh_ges_fb}
Let Assumption \ref{ass.1} hold with $m>0$ and let $\mu = 1/(2L)$. Then, the proximal gradient flow~\eqref{eq.fb_dyn} and DR splitting dynamics~\eqref{eq.dr_dyn} are globally exponentially stable with rate $\rho=\alpha m$.
\end{lemma}

	\vspace*{-1ex}
\section{Main results}\label{sec.main}
		\vspace*{-1ex}

In this section, we examine convergence properties of continuous-time variants of accelerated proximal gradient and DR splitting algorithms for nonsmooth composite optimization problem~\eqref{eq.composite}. In Theorems~\ref{theorem.quad_cvx} and~\ref{theorem.quad_str}, we show that the corresponding dynamical systems provide acceleration when the smooth part of the objective function $f$ is either convex or strongly convex quadratic function, respectively. Furthermore, in Theorem~\ref{theorem.gen_str}, we demonstrate that properties of accelerated proximal gradient flow dynamics carry over from quadratic to general strongly convex problems. To improve readability, we relegate the proofs of theorems to Sections~\ref{sec.quad}.

Inspired by~\eqref{eq.inertial}, which provides accelerated convergence for smooth problems, we propose the following differential equation to accelerate proximal gradient flow dynamics~\eqref{eq.fb_dyn},
\bseq\label{eq.inertial_2}
\begin{IEEEeqnarray}{rcl}
\xdd \,+\, \gamma\xd \,+\, \alpha G_\mu(x \,+\, \beta\xd)  &~=~& 0. 
\label{eq.acc_fb}
\end{IEEEeqnarray}
Similarly, we introduce
\begin{IEEEeqnarray}{rcl}
	\label{eq.acc_dr}
\zdd \,+\, \gamma\zd \,+\, \alpha G_\mu(\prox_{\mu f}(z \,+\, \beta\zd))  &~=~& 0 
	\\
	\non
x &~=~& \prox_{\mu f}(z)	
\end{IEEEeqnarray}
\eseq
as an accelerated variant of DR splitting dynamics~\eqref{eq.dr_dyn}, where $G_\mu$ is the generalized gradient map associated with problem~\eqref{eq.composite}; see~\eqref{eq.gen_grad} for definition. 

Our analysis of~\eqref{eq.acc_fb} and~\eqref{eq.acc_dr} utilizes Lyapunov functions obtained by replacing $f$ in~\eqref{eq.lyap0} with smooth envelope functions $\cF_{\mu}$ or $\cD_{\mu}$, respectively. This choice is motivated by the variable-metric gradient method interpretations of FB and DR splittings. The lack of convexity of the FB and DR envelopes when $f$ in~\eqref{eq.composite} is not a convex quadratic function is the main challenge in convergence analysis. When $f$ is either a convex or a strongly convex quadratic function, we show that~\eqref{eq.acc_fb} and~\eqref{eq.acc_dr} converge to a solution of problem~\eqref{eq.composite} at an accelerated rate. We then remove this quadratic restriction on $f$ and show that acceleration carries over to general strongly convex problems for FB splitting dynamics~\eqref{eq.acc_fb}.

	\vspace*{-1ex}
\subsection{Convex quadratic $f$}
	\vspace*{-1ex}
	
We first restrict our attention to a nonsmooth optimization problem~\eqref{eq.composite} with quadratic $f$. In Theorem~\ref{theorem.quad_cvx}, we establish accelerated sublinear convergence for convex problems and, in Theorem~\ref{theorem.quad_str}, we show that the convergence rate can be improved to exponential for strongly convex problems. This extends results from the literature from a discrete-time to a continuous-time setting.  

\begin{assumption}\label{ass.quad}
In~\eqref{eq.composite}, $f$ is a convex quadratic function with a strong convexity constant $m \geq 0$ and a Lipschitz continuity constant $L$, and $g$ is a proper closed convex function.
\end{assumption}

\begin{theorem}\label{theorem.quad_cvx}
Let Assumption~\ref{ass.quad} hold and let $\mu\in(0, 1/L)$. If the algorithmic parameters are given by
$(\alpha > 0, \gamma(t)  = 3/(t + 3), \beta(t)  = 1  - \gamma(t))$
then the solutions to~\eqref{eq.acc_fb} and~\eqref{eq.acc_dr} for $t \geq 0$ satisfy
\beq\non
F(p_\mu(x(t)))   -  F(\xs)   \,\leq\,  \tfrac{c_1}{(t \,+\, 3)^2}( \norm{x(0)  -  \xs}_2^2  \,+\,  \norm{\xd(0)}_2^2)
\eeq
where $p_\mu(x) \DefinedAs \prox_{\mu g}(x - \mu\nabla f(x))$, $\xs \in \argmin F(x)$, and $c_1 > 0$ is a constant that depends on the problem parameters.
\end{theorem}

Theorem~\ref{theorem.quad_cvx} establishes accelerated sublinear convergence rate for convex problems with quadratic smooth part of the objective function. This matches the optimal rate $O(1/k^2)$ of discrete-time gradient-based methods~\cite{nes83}. Furthermore, in Theorem~\ref{theorem.quad_str} we show that the use of constant algorithmic parameters for strongly convex problems with quadratic $f$ and known $m > 0$ leads to an accelerated exponential convergence rate. In proof of Theorem~\ref{theorem.quad_str}, we utilize Lemma~\ref{lemma.cvx_env} which provides estimates of the smoothness and strong convexity constants of the FB and DR envelopes under Assumption~\ref{ass.quad}.

\begin{lemma}[Thm.~2.3~\cite{patstebem14} and Proposition 4.6~\cite{gisfae18}]\label{lemma.cvx_env}
Let Assumption~\ref{ass.quad} hold and let $\mu\in(0,1/L)$. The FB and DR envelopes $\cF_{\mu}$ and $\cD_{\mu}$ are convex functions with the following smoothness and strong convexity constants 
\beq\non
\ba{rcl}
\tilde{L} &\asp{=}& 
\left\{
\ba{ll}
2(1 \,-\, \mu m)/\mu & ~~\mbox{\rm for}~{\cF_\mu}  
\\[0.1cm]
(1 \,-\, \mu m)/ \big( \mu(1  \,+\, \mu m)^2 \big) & ~~\mbox{\rm for}~{\cD_\mu}  
\ea
\right.
	\\[0.5cm]
\tilde{m} &\asp{=}& 
\left\{
\ba{ll}
\min\left\{(1 \,-\, \mu m)m, (1 \,-\, \mu L)L\right\}& ~~\mbox{\rm for}~{\cF_\mu} 
\\[0.1cm]
\min\left\{\dfrac{(1 \,-\, \mu m)m}{(1 \,+\, \mu m)^2}, \dfrac{(1 \,-\, \mu L)L}{(1 \,+\, \mu L)^2}\right\} & ~~\mbox{\rm for}~{\cD_\mu}  
\ea
\right.	   
\ea
\eeq 
\end{lemma}

If the condition number $\kappa = L/m$ of $f$ satisfies $\kappa\geq 2$, then Lemma~\ref{lemma.cvx_env} shows that for $\mu = 1/(kL)$ with $k \geq 2$, the condition numbers of the FB and DR envelopes, $\tilde{\kappa} \DefinedAs \Lt/\mt$, are both upper bounded by $2k\kappa$.

\begin{theorem}\label{theorem.quad_str}
Let Assumption~\ref{ass.quad} hold with $m>0$, let $\mu\in(0, 1/L)$, and let $\Lt$ and $\mt$ be the smoothness and strong convexity constants of either the FB or the DR envelope. If the algorithmic parameters in~\eqref{eq.inertial_2} are given by
	$
	(
	\alpha > 0, \gamma = 2\sqrt{\alpha\mt}\big/(\sqrt{\alpha\mt}+1), \beta  = 1 - \gamma
	)
	$
then the solutions to~\eqref{eq.acc_fb} and~\eqref{eq.acc_dr} for $t \geq 0$ satisfy
\beq\non
\norm{x(t) \,-\, \xs}_2^2  \;\leq\; c_2(\norm{x(0) \,-\, \xs}_2^2 \,+\, \norm{\xd(0)}_2^2) \, \rme^{-\rho t}
\eeq
where $\rho  =  \sqrt{\alpha\mt} - \alpha\mt/2$ is the convergence rate and $c_2 > 0$ is a constant that depends on the problem parameters.
\end{theorem}

	\vspace*{-1ex}
\subsection{Strongly convex $f$}
	\vspace*{-1ex}

We next show that an additional restriction on~$\mu$ allows us to obtain an accelerated exponential convergence rate for~\eqref{eq.acc_fb} even when the smooth part of the objective function $f$ in~\eqref{eq.composite} is not quadratic but strongly convex.

\begin{theorem}\label{theorem.gen_str}
Let Assumption~\ref{ass.1} hold and let $\mu\in\left( 0,\sqrt{\gamma\beta}/(2L)\right]$. If the algorithmic parameters in~\eqref{eq.acc_fb} are given by
	$
	(
	\alpha > 0, \gamma = 2\sqrt{\alpha m}\big/(\sqrt{\alpha m}  + 1), \beta = 1 - \gamma
	)
	$
then the solutions to~\eqref{eq.acc_fb} for $t \geq 0$ satisfy
\beq\non
\norm{x(t) \,-\, \xs}_2^2  \;\leq\; c_3(\norm{x(0) \,-\, \xs}_2^2 \,+\, \norm{\xd(0)}_2^2)\rme^{-\rho t}
\eeq
where $\rho  =  \sqrt{\alpha m} - \alpha m/2$ is the convergence rate, and  $c_3 > 0$ is a constant that depends on the problem parameters. 
\end{theorem}
Setting $\alpha = 1/\Lt$ and $\alpha = 1/L$ in Theorems~\ref{theorem.quad_str} and~\ref{theorem.gen_str}, respectively, gives an exponential convergence rate, $\rho\sim O(1/\sqrt{\kappa})$. This rate matches the best-achievable rate of discrete-time gradient-based methods for strongly convex problems~\cite{nes03}. Moreover, for $\alpha = 1/L$, the upper bound on the penalty parameter~$\mu$ in Theorem~\ref{theorem.gen_str} is of $O(1/(L\sqrt[4]{\kappa}))$ and is smaller than the upper bound on $\mu$ in Theorem~\ref{theorem.quad_str} by a factor of $O(1/\sqrt[4]{\kappa})$.

\begin{remark}
For the rate-preserving discretization proposed in~\cite{muejor19} to be applicable to accelerated dynamics~\eqref{eq.acc_fb}~and~\eqref{eq.acc_dr}, we require $\beta + \gamma = 1$ in Theorems~\ref{theorem.quad_cvx}-\ref{theorem.gen_str}. For general convex problems, this additional constraint on algorithmic parameters together with the lack of convexity of envelope functions makes the analysis challenging. For~\eqref{eq.acc_fb}, we obtain acceleration by introducing additional constraints on  the penalty parameter $\mu$ and two novel inequalities that substitute the strong convexity inequalities; see Lemma~\ref{lemma.ineq} for details.
\end{remark}

	\vspace*{-2ex}
\section{Proof of main results}\label{sec.quad}
	\vspace*{-1ex}
	
We next present a system-theoretic framework that can be specialized to particular problem setups addressed in Theorems~\ref{theorem.quad_cvx}-\ref{theorem.gen_str}. To simplify notation, we assume without loss of generality that in~\eqref{eq.acc_fb}, $\min_x \cFm(x) = 0$ and that $\xs = 0$ is a minimizer; similarly, for \eqref{eq.acc_dr}, $\min_z \cDm(z) = 0$ and $\zs = 0$ is a minimizer.

\textbf{System representation.} We rewrite~\eqref{eq.inertial_2} as a feedback of a linear system with a static nonlinear block,
\bseq\label{eq.ss-fbk}
\beq\label{eq.dyn_sys}
\dot{\pp}   \,=\,   A\pp \,+\, Bu,
\quad
y \,=\,  C\pp,
\quad
u  \,=\,  \Delta(y)
\eeq
where $\pp = [~x^T~\xd^T~]^T$ and $\pp = [~z^T~\zd^T~]^T$ denote the respective state vectors for~\eqref{eq.acc_fb} and~\eqref{eq.acc_dr}. On the other hand, the system matrices are given by
	\beq\label{eq.sys_matrices}
\mleft[
\begin{array}{c|c}
  A & B \\
  \hline
  C & 0
\end{array}
\mright]
   \,=\,   
 \mleft[
 \begin{array}{cc|c}
   0 & I & 0 \\
   0 & -\gamma I & -\alpha I
   \\
   \hline
   I & \beta I & 0
 \end{array}
 \mright].
\eeq
Here, the partition is done conformably with the partition of the state vector $\pp = [~\pp_1^T~\pp_2^T~]^T$ and the corresponding nonlinear terms are
	\beq\label{eq.sys_nonlinear}
	\Delta(y) 
	\, = \,
	\left\{
	\ba{ll}
	G_\mu(y), & \text{FB splitting}
	\\[-0.05cm]
	G_\mu(\prox_{\mu f}(y)), 
	& 
	\text{DR splitting.}
	\ea
	\right.
	\eeq
	\eseq
	
	\vspace*{-2ex}
\subsection{Theorems~\ref{theorem.quad_cvx} and~\ref{theorem.quad_str}: convex quadratic $f$}
\vspace*{-1ex}

We first prove our results for convex quadratic $f$ in~\eqref{eq.composite}, $f(x) = (1/2)  x^TQx + q^Tx$, where $Q\in\reals{n\times n}$ is a positive semidefinite matrix and $q \in \reals{n}$. In this case, the Lipschitz continuity and strong convexity parameters $L$ and $m$ are given by the largest and smallest eigenvalues of $Q$. 

\textbf{Lyapunov function.}
Inspired by~\eqref{eq.lyap0} which was used to establish accelerated convergence rates for~\eqref{eq.inertial} in~\cite{muejor19}, we propose the Lyapunov function candidate for~\eqref{eq.inertial_2},
\beq\label{eq.lyap}
V(t; \pp) \, = \, \alpha\cE(\pp_1) \, + \, (1/2)\norm{\theta(t)\pp_1  +  \pp_2}_H^2
\eeq
where the envelope $\cE$ (either $\cF_\mu$ or $\cD_\mu$) and the associated weight matrix $H$ are determined by the selection of the nonlinear term $\Delta$ in~\eqref{eq.sys_nonlinear}, and $\norm{x}_H^2 \DefinedAs \inner{x}{x}_H = x^THx$. Our aim is to show that $\dot{V}(t) + \theta (t) V(t)\leq 0$ for $t\geq t_0$ which in conjunction with the Gr{\"o}nwall inequality implies $V(t)\leq V(t_0)\exp(-\int_{t_0}^{t}\theta(\tau)\rmd\tau)$. The inequalities presented in the theorems are then obtained by leveraging various properties of associated envelopes and the proximal operators. Hence, $\theta (t)$ in~\eqref{eq.lyap} determines the convergence rate. As we elaborate in the forthcoming subsections, in Theorem~\ref{theorem.quad_cvx} we set $\theta(t) = 2/(t+3)$ and in Theorem~\ref{theorem.quad_str} we set $\theta(t) = \rho$ to establish sublinear and exponential convergence rates, respectively.  

\textbf{Characterization of FB and DR envelopes.}
To upper bound $\theta V$, we first derive an upper bound on the envelope function $\cE$. For a (strongly) convex quadratic $f$, the FB and DR envelopes $\cE$ are also (strongly) convex; see Lemma~\ref{lemma.cvx_env}. Thus, for any $x,\xp\in\reals{n}$, 
\beq\label{eq.str_ineq_eps}
\cE(x) \, - \, \cE(\xp)   \, \leq \, \inner{\nabla \cE(x)}{x \,-\, \xp} \, - \, \dfrac{\mt}{2} \, \norm{x \,- \, \xp}_2^2
\eeq 
where  $\mt$ is the parameter of strong convexity of $\cE$. Moreover, \eqref{eq.fb_grad} and~\eqref{eq.dr_grad2} yield $\nabla \cE(x) = H\Delta(x)$. For the FB envelope $\cE = \cF_{\mu}$, $\Delta(x) = G_\mu(x)$, and
\beq\non
H  \, = \,  I \, - \, \mu\nabla^2f(x) \, = \, I  \, - \,  \mu Q.
\eeq
On the other hand, for the DR envelope $\cE = \cD_{\mu}$, $\Delta(x) = G_\mu(\prox_{\mu f}(x))$ and the matrix $H$ is determined by 
\beq\non
H  \,=\,  2\nabla \prox_{\mu f}(x)  \,-\,  I  \,=\,  (I  \,-\,  \mu Q)(I  \,+\,  \mu Q)^{-1}.
\eeq
In both cases, $0 \prec H\preceq I$ for $\mu\in(0,1/L)$, which together with~\eqref{eq.str_ineq_eps}  implies that for any $x,\xp\in\reals{n}$
\beq\label{eq.delta_char}
\cE(x) \, - \, \cE(\xp)   \, \leq \, \inner{\Delta(x)}{x \,-\, \xp}_H \, - \, \dfrac{\mt}{2} \, \norm{x \,- \, \xp}_H^2.
\eeq 

\textbf{An upper bound on the Lyapunov function.}
We first rewrite the Lyapunov function candidate using system representation~\eqref{eq.dyn_sys} as
\beq\label{eq.lyap2}
V(t; \pp)  \,=\,  \alpha\cE(C_1\pp) \,+\, (1/2)\norm{R(t) \pp}_H^2
\eeq
where $C_1 = [I~0]$, $R(t) = [\theta(t)I~I]$. We define $C_2  = [0~\beta I]$, $v = \Delta(C_1\pp)$, and recall our assumption that $\pp_1^\star = C_1 \pp^\star = 0$ and $\cE(\pp_1^\star) = 0$. Using inequality~\eqref{eq.delta_char} twice, for $(x,\xp) = (C\pp, C_1 \pps)$ and $(x,\xp) = (C_1\pp, C\pp)$, along with $u = \Delta (C \psi)$ and $C = C_1 + C_2$ yield,
\beq\non
\ba{rcl}
\cE(C\pp) &\asp{\leq}& \inner{u}{C\pp}_H   \,-\, (\mt/2)\norm{C\pp}_H^2
\\ 
\cE(C_1\pp)   \,-\,  \cE(C\pp) &\asp{\leq}& -\inner{v}{C_2\pp}_H   \,-\,  (\mt/2)\norm{C_2\pp}_H^2.
\ea
\eeq
Finally, using the definition of $V$ in~\eqref{eq.lyap2} and summing the last two inequalities gives an upper bound on $\theta V$,
\beq\label{eq.upper_lyap}
\ba{rcl}
\!\!\!\!\!\theta V  &\asp{\leq}& \alpha\theta \left( \inner{u}{C\pp}_H   \,-\,  \inner{v}{C_2\pp}_H\right) 
\\ 
&\asp{}&-\,  \dfrac{\alpha\mt\theta}{2}\left( \norm{C\pp}_H^2  \,+\,  \norm{C_2\pp}_H^2\right) \,+\, \dfrac{ \theta}{2}\norm{R\pp}_H^2 .
\ea
\eeq
We next examine the derivative of $V$ along the solutions of~\eqref{eq.ss-fbk} and use~\eqref{eq.upper_lyap} to obtain an upper bound on $\dot{V}+\theta V$.

\textbf{Lyapunov analysis.}
The derivative of~\eqref{eq.lyap2} along the solutions of~\eqref{eq.ss-fbk} is determined by
\beq\label{eq.lyap_dot}
\dot{V} = \alpha\inner{v}{C_1(A\pp + Bu)}_H + \inner{R\pp}{\dot{R}\pp + R(A\pp + Bu)}_H
\eeq 
where $\dot{R} = [(\rmd\theta (t)/\rmd t) I~0]$. Finally, combining~\eqref{eq.lyap_dot} with~\eqref{eq.upper_lyap} leads to
\beq\label{eq.master}
\dot{V}  \,+\,  \theta V  \,\leq\,  \frac{1}{2}\thbo{\pp \\[-0.75cm]}{u \\[-0.75cm]}{v}^T\thbth{\Pi}{\Gamma}{\Lambda \\[-0.75cm]}{\Gamma^T}{0}{0 \\[-0.75cm]}{\Lambda^T}{0}{0}\thbo{\pp \\[-0.75cm]}{u \\[-0.75cm]}{v}
\eeq
where matrices $\Pi$, $\Gamma$, and $\Lambda$ differ for convex and strongly convex cases. Now, we prove that the matrix on the right-hand-side of~\eqref{eq.master} is negative semidefinite for the parameters given in Theorems~\ref{theorem.quad_cvx} and~\ref{theorem.quad_str}. 

	\vspace*{-1ex}
\subsubsection{Theorem 1: convex quadratic $f$} \label{sec.nonstr_quad}
	\vspace*{-1ex}
	
For a (non-strongly) convex quadratic function $f$, $\mt=0$, the upper bound on $\theta V$ given in~\eqref{eq.upper_lyap} simplifies to
\beq\non
\theta V  \,\leq\,  \alpha\theta \left( \inner{u}{C\pp}_H   \,+\,  \inner{v}{-C_2\pp}_H\right)  \,+\, (\theta/2)\norm{R\pp}_H^2.
\eeq 
We let $P \DefinedAs R^THR$ and utilize this upper bound and~\eqref{eq.lyap_dot} to write $\Pi$ in~\eqref{eq.master} as 
\beq\label{eq.pi}
\ba{rcl}
\Pi&\asp{=}& 
A^TP  \,+\,  PA  \,+\,  \theta P \,+\,  \dot{P} 
\\ 
&\asp{=}& \tbt{2\dot{\theta}\theta  \,+\, \theta^3}{\dot{\theta} \,+\, 2\theta^2 \,-\, \theta\gamma \\[-0.65cm]}{\dot{\theta} \,+\, 2\theta^2 \,-\, \theta\gamma}{3\theta \,-\,  2\gamma} \otimes H
\ea
\eeq
where $\dot{P} = \dot{R}^THR + R^TH\dot{R}$ and $\otimes$ is the Kronecker product. A necessary condition for negative semidefiniteness of the matrix on the right-hand-side of~\eqref{eq.master} is the negative semidefiniteness of~$\Pi$. One way to simplify the conditions that guarantee $\Pi\preceq0$ is to set the off-diagonal blocks in~\eqref{eq.pi} to zero, i.e.,
\beq\label{eq.cross1}
\dot{\theta}  \,+\, 2\theta^2  \,-\,  \theta\gamma  \,=\,  0.
\eeq
When~\eqref{eq.cross1} is satisfied, $\Pi\preceq 0$ if its diagonal blocks are negative semidefinite. The first diagonal block requires $\theta \leq 2/(t + r)$ and the largest~$\theta$ that satisfies this condition is $\theta (t) = 2/(t+r)$ for some $r>0$. Substituting this value to~\eqref{eq.cross1} yields $\gamma (t) = 3/(t+r)$. For this choice of parameters, $\Pi \equiv 0$. While any choice of $r > 0$ works here, different constants $c_1$ in Theorem~\ref{theorem.quad_cvx} are obtained for different values of $r$. In what follows, we set $r=3$ to ensure positivity of $\beta(t)$ in Theorem~\ref{theorem.quad_cvx} for any $t\geq 0$. 

Since the off-diagonal blocks in~\eqref{eq.master} are given by
\beq\label{eq.cross}
\ba{rclll}
\Gamma &\asp{=}& PB  \,+\,  \alpha\theta C^TH &\asp{=}&\tbo{0}{-\alpha(1  \,-\,  \theta\beta)H}
\\ 
\Lambda &\asp{=}& \dfrac{\alpha}{\beta}(1 \,-\,  \theta\beta) C_2^TH &\asp{=}& \tbo{0}{\alpha(1 \,-\, \theta\beta)H}  \,=\, -\Gamma.
\ea
\eeq
for the above selection of $\theta$ and $\gamma$,~\eqref{eq.master} simplifies to
\beq\label{eq.master_1}
\ba{rcl}
\dot{V} \,+\,  \theta V  &\asp{\leq}& -(\alpha/\beta)(1 \,-\, \theta\beta)\inner{u \,-\, v}{C_2\pp}_H .
\ea
\eeq
Moreover, monotonicity of the gradient of the envelope function implies $\inner{u - v}{C_2\pp}_H\geq 0$. Thus, the right-hand-side of~\eqref{eq.master_1} is upper bounded by zero and the use of Gr{\"o}nwall inequality yields
$V(t) \leq V(0)/(t + 3)^2$, $t\geq0$.
The inequality in Theorem~\ref{theorem.quad_cvx} is obtained as follows. We first use the smoothness of the envelope function together with our assumption that $\min_x\cF_\mu(x) = 0$ and $\xs \in \argmin_x\cF_\mu(x)$ to obtain
\beq\non
\cE(\psi_1(0))  \,\leq\, \tfrac{\Lt}{2} \norm{\psi_1(0)}^2  \,\leq\, \tfrac{\Lt}{2} \norm{\psi(0)}^2
\eeq
where $\Lt$ is given in Lemma~2.  Then, application of the triangle inequality, $H\preceq I$, and the fact that $\theta(t)\leq 1$ for $t\geq 0$ yield the upper bound on the Lyapunov function,
\begin{align}
	V(\psi(0))  &\,=\, \tfrac{\alpha \Lt}{2} \norm{\psi(0)}^2 \,+\,  \tfrac{1}{2}\norm{\theta(t)\pp_1(0) \,+\, \pp_2(0)}_H^2
	\non
	\\[\nvs]
	\label{eq.upp}
	&\,\leq\, 
	(\tfrac{\alpha \Lt}{2}  \,+\, 1)\norm{\psi(0)}^2.
\end{align}
We note that for both~\eqref{eq.acc_fb} and~\eqref{eq.acc_dr}, $\cE(\pp_1) = \cF_\mu(x)$. Hence, a lower bound on the Lyapunov function can be obtained using the identity~\cite[Theorem~2.2]{patstebem14}, 
$F(p_\mu(x))\leq \cF_{\mu}(x)$ for any $x \in \reals{n}$. This lower bound in conjunction with~\eqref{eq.upp} leads to the inequality in Theorem~1 for~\eqref{eq.acc_fb} with constant $c_1 =1/\Lt + 1/\alpha$.

In the case of~\eqref{eq.acc_dr} for which $\psi = [z^T~\dot{z}^T]^T$, while the lower bound on the Lyapunov function stays the same, the upper bound~\eqref{eq.upp} can be expressed in terms of $x$ and $\dot{x}$ using the following relations~\cite{patstebem14a},
\beq\non
\ba{rcl}
x &\asp{=}&  \prox_{\mu f}(z)  \,=\,  (I \,+\, \mu Q)^{-1}z,
\\[\nvs] 
\xd  &\asp{=}&  \nabla \prox_{\mu f}(z)\dot{z}  \,=\, (I \,+\, \mu Q)^{-1}\dot{z}
\ea
\eeq
which lead to the inequality in Theorem~1 with constant $c_1 = (1 + \mu L)(1/\Lt + 1/\alpha)$.

	\vspace*{-1ex}
\subsubsection{Theorem 2: strongly convex quadratic $f$}\label{sec.quad_str}
	\vspace*{-1ex}
	
For a strongly convex $f$, the use of constant $\alpha$, $\beta$, and $\gamma$ in~\eqref{eq.inertial_2} leads to time-invariant dynamics~\eqref{eq.dyn_sys}. In this case, the Lyapunov function candidate does not explicitly depend on time, i.e., $R$ in~\eqref{eq.lyap2} becomes a constant matrix. We again let $P = R^THR$ and utilize~\eqref{eq.upper_lyap} and~\eqref{eq.lyap_dot} to write $\Pi$ in~\eqref{eq.master} as
\beq\non
\ba{rcl}
\Pi&\asp{=}& A^TP \,+\, PA \,+\, \theta P  \,-\,  \alpha\mt\theta(C^THC \,+\, C_2^THC_2).
\ea
\eeq
Since the off-diagonal blocks $\Gamma$ and $\Lambda$ in~\eqref{eq.master} remain the same as in~\eqref{eq.cross}, we have
\beq\label{eq.off_diag}
\inner{\pp}{\Gamma u   \,+\,  \Lambda v} 
 \,=\,   
  -(\alpha/\beta)(1 \,-\, \theta\beta)\inner{C_2\pp}{u-v}_H.
\eeq
Moreover, strong monotonicity of gradient of the envelope function along with $H\preceq I$ implies $\inner{u  -  v}{C_2\pp}_H \geq \mt\norm{C_2\pp}_H^2$ which together with~\eqref{eq.off_diag} simplifies~\eqref{eq.master} to
\beq\label{eq.master2}
\dot{V}  \,+\,  \theta V  \,\leq\, (1/2)\pp^T(\Pi \,-\,  \Upsilon)\pp 
\eeq
where
\beq\non
\Upsilon  \,=\,  (2\alpha\mt/\beta)(1 \,-\,  \theta\beta) C_2^THC_2.
\eeq
Substitution of the system matrices in~\eqref{eq.sys_matrices} gives
\beq\non
\Pi   -  \Upsilon
= 
\tbt{ \theta(\theta^2  -  \alpha\mt)}{\theta(2\theta  -  \gamma  -  \alpha\mt\beta)}{\theta(2\theta  -  \gamma  -  \alpha\mt\beta)}{3\theta  -  2\gamma  -  2\alpha\mt\beta }\otimes H. 
\eeq
The conditions that guarantee $\Pi-\Upsilon\preceq 0$ can be simplified by setting off-diagonal blocks in $\Pi-\Upsilon$ to zero, i.e.,
\beq\label{eq.cross21}
\theta   \,=\,   (\gamma \,+\, \alpha\mt\beta)/2.
\eeq
When~\eqref{eq.cross21} is satisfied, the second diagonal block in $\Pi  - \Upsilon$ is clearly negative definite. Since $\beta = 1-\gamma$, the negative semidefiniteness of the first diagonal block requires
$
\gamma \leq(2\sqrt{\alpha\mt}  -  \alpha\mt)/(1 - \alpha\mt)
$
which is trivially satisfied by the choice of $\gamma= (2\sqrt{\alpha\mt} - 2\alpha\mt)/(1 - \alpha\mt)$ given in the theorem. Hence, the right-hand-side of~\eqref{eq.master2} is negative semidefinite, leading to $V(t) \leq V(0)\rme^{-\theta t}$ for $ t \geq 0$.

The strong convexity and smoothness of the envelope functions yield quadratic lower and upper bounds on the Lyapunov function, which completes the proof for~\eqref{eq.acc_fb}. For~\eqref{eq.acc_dr},  the proof is completed  as in Section~\ref{sec.nonstr_quad}, i.e., utilizing linearity of the proximal operator together with $\xd = \nabla\prox_{\mu f}(z)\zd$ and~\eqref{eq.prox_grad}.

	\vspace*{-1ex}
\subsection{Theorem 3: strongly convex $f$}\label{sec.gen_str}
	\vspace*{-1ex}

To prove our results for general strongly convex~$f$, we use system representation given by~\eqref{eq.ss-fbk} with constant parameters but with a different Lyapunov function and a different characterization of the FB envelope.

\textbf{Lyapunov function.} For general strongly convex $f$, since the weight matrix $H (t)  =  I - \mu\nabla^2f(x (t))$ is time-dependent, in lieu of~\eqref{eq.lyap2} we propose the following time-invariant Lyapunov function candidate
\beq\label{eq.lyap3}
V(\pp) \,=\, \alpha \cF_\mu(C\pp) \,+\, (1/2)\pp^TP\pp
\eeq 
where $C$ is given in~\eqref{eq.sys_matrices}, $P = R^TR$, and $R$ is given in~\eqref{eq.lyap2}. In contrast to~\eqref{eq.lyap2}, the envelope function in~\eqref{eq.lyap3} is evaluated at $y =C\pp$ and the quadratic form does not use the time-dependent weight matrix~$H$. 

\textbf{Characterization of the FB envelope.}
For general strongly convex problems, the FB envelope is not necessarily convex and the upper bounds~\eqref{eq.str_ineq_eps} and~\eqref{eq.delta_char} no longer hold. In Lemmas~\ref{lemma.ineq}, we derive two novel inequalities that substitute~\eqref{eq.str_ineq_eps}; see Appendix for the proof.

\begin{lemma}\label{lemma.ineq}
Let Assumption~\ref{ass.1} hold and let $\mu\in(0,1/L)$. For every $x, \xp \in \reals{n}$ and $\xs \in \argmin F(x)$,  FB envelope satisfies   
\bseq\label{eq.lemma}
\begin{gather}
\label{eq.lemma.1}
\ba{l}
\!\!\!\!\!\cF_\mu(x) \,-\, F(\xp)\,\leq\, 
\\[-0.05cm] 
\!\!\!\!\!\inner{G_\mu(x)}{x \,-\, \xp}  \,-\, \frac{m}{2}\norm{x \,-\, \xp}_2^2 \,-\, \frac{\mu}{2}\norm{G_\mu(x)}_2^2.
\ea
\\[0.cm]
\label{eq.lemma.2}
\cF_\mu(x) \,-\, F(\xs)  \,\geq\, \tfrac{m^2(1 \,-\, \mu L)}{2L} \norm{x \,-\, \xs}_2^2.
\end{gather}
\eseq
\end{lemma}
We note that in contrast to strong convexity inequality~\eqref{eq.str_ineq_eps}, in Lemma~\ref{lemma.ineq}, the gradient of $\cF_\mu$ is not involved and we characterize the difference between envelope $\cF_\mu$ and objective function $F$. Still, Lyapunov function~\eqref{eq.lyap3} can be upper bounded via Lemma~\ref{lemma.cvx_env} based on the equivalence given in~\eqref{eq.env_equivalence}.

\textbf{An upper bound on the Lyapunov function.}
The following upper bound on~\eqref{eq.lyap3} can be directly obtained using Lemma~\ref{lemma.ineq},~\eqref{eq.env_equivalence}, and system representation~\eqref{eq.dyn_sys},
\beq\label{eq.lyap_upper2}
\ba{rcl}
\theta V &\asp{\leq}& \alpha\theta\inner{u}{C\pp}    -  \frac{\alpha m\theta}{2}\norm{C\pp}_2^2  -  \frac{\alpha\mu\theta}{2}\norm{u}_2^2  +  \tfrac{\theta}{2}\pp^TP\pp.
\ea
\eeq

\textbf{Lyapunov analysis.}
The derivative of~\eqref{eq.lyap3} along the solutions of~\eqref{eq.acc_fb} is determined by
\beq\non\label{eq.lyap_dot2}
\dot{V}  \,=\, \alpha\inner{H u}{CA\pp \,+\, CBu} \,+\, \inner{P\pp}{A\pp \,+\, Bu}.
\eeq
Summing $\dot{V}$ with the upper bound given in~\eqref{eq.lyap_upper2} yields
\beq\non
\dot{V} \,+\, \theta V \,\leq\, \frac{1}{2}\tbo{\pp}{u}^T\tbt{\Pi}{\Gamma + \Lambda H}{\Gamma^T + H \Lambda^T}{-2\alpha^2\beta H}\tbo{\pp}{u}
\eeq
where
\beq\label{eq.matrices}
\ba{rcl}
\Pi &\asp{=}& A^TP \,+\, PA \,+\,  \theta P \,-\, \alpha m\theta C^TC 
\\ 
\Gamma &\asp{=}& PB \,+\, \alpha\theta C^T
\\ 
\Lambda &\asp{=}&  \alpha A^TC^T.
\ea
\eeq
Since $H \succ 0$, the Schur complement implies that $\dot{V} + \theta V \leq 0$ if and only if 
\beq\label{eq.gen_ineq}
2\alpha^2\beta\Pi + (\Gamma\Lambda^T + \Lambda\Gamma^T) + \Gamma H^{-1}\Gamma^T + \Lambda H \Lambda^T \preceq 0.
\eeq
Substitution of system matrices~\eqref{eq.sys_matrices}  into~\eqref{eq.matrices} gives
\bseq
\begin{align}
&\Gamma\Lambda^T \,+\, \Lambda\Gamma^T  \,=\, \tbt{0}{0}{0}{-2\alpha^2(1  \,-\,  \theta\beta)(1 \,-\, \gamma\beta) I}\label{eq.cross2}
\\ 
&\Gamma H^{-1}\Gamma^T \,+\, \Lambda H \Lambda^T  \,=\, \nonumber
\\ 
&\tbt{0}{~0}{0}{~\alpha^2 \! \left( (1  -  \theta\beta)^2H^{-1} + (1 -\ \gamma\beta)^2H \right)}.\label{eq.cross3}
\end{align}
\eseq
The eigenvalue decomposition of $H$ can be used to show
\bseq
\begin{gather}
\label{eq.square}
\Gamma H^{-1}\Gamma^T \,+\, \Lambda H \Lambda^T 
\,\preceq\,
\tbt{0}{0}{0}{\alpha^2\eta I}
\\
\label{eq.eta}
\!\! \eta  \DefinedAs \maximize\limits_{\sigma\,\in\,[m,L]}~\dfrac{(1 - \theta\beta)^2 }{1  -  \mu\sigma} \,+\,  (1 -  \gamma\beta)^2(1 - \mu\sigma).
\end{gather}
\eseq
If $(1-\theta\beta)\geq(1-\mu\sigma)(1-\gamma\beta)$ for all $\sigma\in[m,L]$, then the objective function in~\eqref{eq.eta} is increasing for $\sigma\in[m,L]$ and the maximum occurs at $\sigma = L$,
\beq\label{eq.equiv_prob}
\eta  \,=\, 
(1 \,-\, \theta\beta)^2/(1 \,-\, \mu L) \,+\, (1 \,-\, \gamma\beta)^2(1 \,-\, \mu L).
\eeq
Combining~\eqref{eq.cross2} with~\eqref{eq.square} and~\eqref{eq.equiv_prob} yields
\beq\non
\ba{l}
(\Gamma\Lambda^T \,+\, \Lambda\Gamma^T) \,+\, \Gamma H^{-1}\Gamma^T \,+\, \Lambda H \Lambda^T  \,\preceq\,  
\\ 
  \tbt{0}{0}{0}{\frac{\alpha^2}{1  -  \mu L} \big((1  -  \theta\beta)   -  (1  -  \mu L)(1  -  \gamma\beta)\big)^2 I}.
\ea
\eeq
Moreover, since the entries of the matrix $\Pi$ are given by 
\beq\non
\Pi  \,=\,   \tbt{\theta(\theta^2  \,-\, \alpha m)}{\theta(2\theta  \,-\, \gamma \,-\, \alpha m\beta)}{\theta(2\theta  \,-\, \gamma \,-\, \alpha m\beta)}{3\theta \,-\, 2\gamma \,-\, \alpha m \theta\beta^2}\otimes I
\eeq
setting the off-diagonal blocks in $\Pi$ to zero yields
\beq\label{eq.cond_cross}
\theta  \,=\, (\gamma \,+\, \alpha m\beta)/2.
\eeq
When~\eqref{eq.cond_cross} is satisfied, verifying $\dot{V} + \theta V \leq 0$, or equivalently~\eqref{eq.gen_ineq}, amounts to checking these three conditions,
\begin{enumerate}[label=(\roman*)]
	\item $(1-\theta\beta)\geq(1-\mu\sigma)(1-\gamma\beta)$ for all $\sigma\in[m,L]$
	\item $\theta^2\leq\alpha m$
	\item $\frac{((1  -  \theta\beta)  - (1  -  \mu L)(1  -  \gamma\beta))^2}{2\beta(1  -  \mu L)} 
	\, \leq \, 
	\alpha m\theta \beta^2 + 2\gamma - 3\theta$.
\end{enumerate}
While condition (i) is necessary for upper bound~\eqref{eq.equiv_prob}, conditions~(ii) and~(iii) ensure negative semidefiniteness of the diagonal blocks in~\eqref{eq.gen_ineq}. 

For $\gamma$ given in Theorem~\ref{theorem.gen_str} and~$\beta = 1-\gamma$, we obtain $\theta = \sqrt{\alpha m} - \alpha m/2$ using~\eqref{eq.cond_cross}. Thus, $\theta \leq \gamma$ and conditions~(i)~and~(ii) are trivially satisfied. For large values of $\kappa = L/m$, the numerator of the left-hand-side in condition~(iii) is dominated by $\mu L$. On the other hand, the largest term in the product of the denominator and the right-hand-side is $\gamma\beta$, which suggests that condition~(iii) holds if $\mu\sim O(\sqrt{\gamma\beta}/L)$. Due to space limitations, we use computational tools to verify this hypothesis. With the above choices of $(\gamma,\beta,\theta)$ and $\mu L = \sqrt{\gamma\beta}/2$, condition~(iii) can be written as $wh(w)\leq 0$ where $w= \sqrt{\alpha m} \in[0,1]$ and $h$: $[0,1]\to\reals{}$. In Fig.~\ref{fig.1}, we show that $h$ is negative for all $w\in[0,1]$. Moreover, since~(i) implies that $w h (w)$ is an increasing function of $\mu L$, condition~(iii) is satisfied for all $\mu L\in(0, \sqrt{\gamma\beta}/2]$. 

\begin{figure}[h]
  \centering
  {
    \begin{tabular}{cc}
      \begin{tabular}{c}
        \vspace*{0.5cm}
        \normalsize{\rotatebox{90}{$h(w)$}}
      \end{tabular}
      &
      \hspace*{-0.85cm}
      \begin{tabular}{c}
        \includegraphics[width=.2\textwidth]{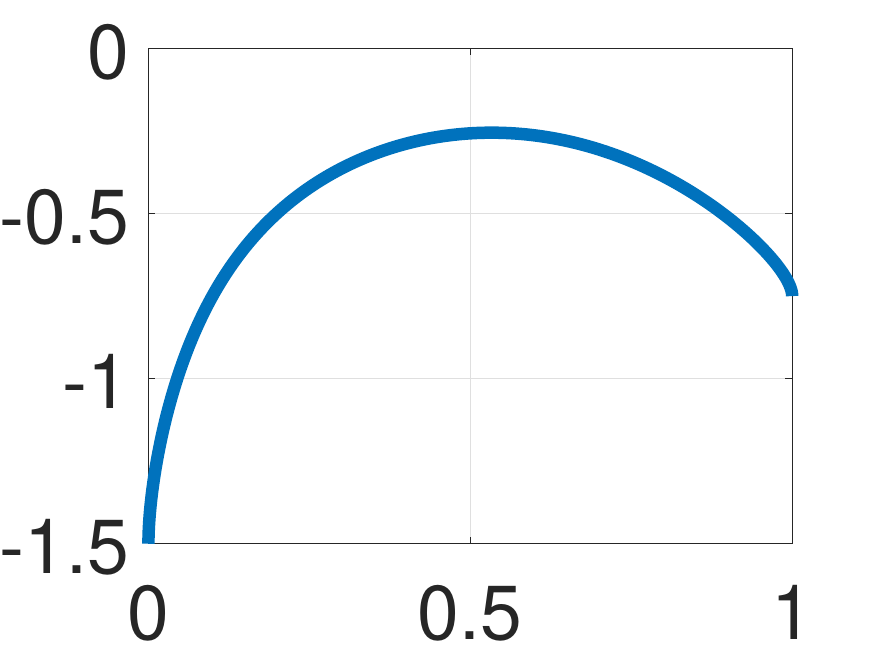}
        \\[-0.1cm]
        $w$
      \end{tabular}   
    \end{tabular}
  }
  \caption{Plot of the function $h$ that demonstrates $h < 0$. Thus, $w h(w)\leq0$ for $w\in[0,1]$ and condition~(iii) is satisfied.}
  \label{fig.h}
\end{figure}
As a result, using Gr{\"o}nwall inequality, we conclude that $V(t) \leq V(0)\rme^{-\theta t}$ for $t \geq 0$ along the solutions of~\eqref{eq.acc_fb}. A quadratic upper bound on the Lyapunov function can be obtain by using the smoothness of the FB envelope as in the proof of Theorem~\ref{theorem.quad_cvx}. Moreover, since $C^TC + P$ is invertible, inequality~\eqref{eq.lemma.2} in Lemma~\ref{lemma.ineq} yields a quadratic lower bound on $V(t)$. These upper and lower quadratic bounds completes the proof.

	\vspace*{-2ex}
\section{Computational experiments}
	\label{sec.simulation}
	
	\vspace*{-1ex}
We use Matlab's function ode45 to evaluate performance of our approach. We choose zero initial conditions, $\alpha=1/L$, and provide comparison with Accelerated ADMM and Tseng's Splitting~\cite{frarobvid21}. Parameters for these discrete-time algorithms are selected as described in~\cite[Sec.~IX.C]{frarobvid21} with decaying (Example~1) and constant (Example~3) damping. Error trajectories of these are plotted at times determined by $t_k = kh/\alpha$, where $k$ is the iteration index, $h=1/L$ is the stepsize of the discrete-time algorithm, and $\alpha$ is the time scale in~\eqref{eq.inertial_2}. While all accelerated algorithms exhibit similar convergence behavior in the convex case (Example~1), our proposed accelerated proximal gradient flow achieves slightly faster convergence in the strongly convex problem (Example~3).

	\vspace*{-1ex}
\subsection{Example 1: $\ell_1$-regularized least squares}
	\vspace*{-1ex}
Let us examine an $\ell_1$-regularized least squares problem,
\beq
\minimize\limits_{x}~(1/2)\norm{Ex  \,-\,  q}_2^2  \,+\, \lambda\norm{x}_1
	\label{eq.ls-ell1}
\eeq
where the problem data $E\in\reals{100\times2000}$, $q \in \reals{100}$, $\lambda>0$ are generated as described in~\cite{odcan15}. Since $E^T E$ is a singular matrix, this problem is not strongly convex. Thus, the time-varying parameters given in Theorem~\ref{theorem.quad_cvx}  with $\mu = 1/(2L)$ are used to illustrate improved sublinear convergence rate of accelerated proximal gradient flow and DR splitting dynamics~\eqref{eq.acc_fb} and~\eqref{eq.acc_dr} in Fig.~\ref{fig.2}.

	\vspace*{-1ex}
\subsection{Example 2: Box-constrained quadratic program}
	\vspace*{-1ex}
	
For a box-constrained quadratic program,
\beq
\label{eq.qp-box}
\ba{l}
\minimize\limits_x~(1/2)x^TQx \,+\, q^Tx  \,+\, I_\cC(x)
\ea
\eeq
where $I_\cC(x)$ is the indicator function of the convex set $\cC = \set{x \, | \, l\leq x\leq u}$, we generate problem data with positive definite $Q\in\reals{500\times500}$, and vectors $q, l,u\in\reals{500}$ as described in~\cite{gonkarros13} and set the condition number $\kappa = 10^5$. Linear convergence of the algorithms for constant parameters given in Theorem~\ref{theorem.quad_str} is illustrated in Fig.~\ref{fig.2}.

	\vspace*{-1ex}
\subsection{Example 3: $\ell_1$-regularized logistic regression}
\begin{figure*}[t]
	\centering
	\begin{tabular}{c@{\hspace{-0.4 cm}}c@{\hspace{-0.3 cm}}c@{\hspace{-0.3 cm}}c@{\hspace{-0.3 cm}}c@{\hspace{-0.3 cm}}c}
		&\subfigure[]{\label{fig.1}}
		&&
		\subfigure[]{\label{fig.2}}
		&&
		\subfigure[]{\label{fig.3}}
		\\[-.2cm]
		\begin{tabular}{c}
			\vspace{.25cm}
			\normalsize{\rotatebox{90}{$F(x(t)) - F(\xs)$}}
		\end{tabular}
		&
		\begin{tabular}{c}
			\includegraphics[width=0.28\textwidth]{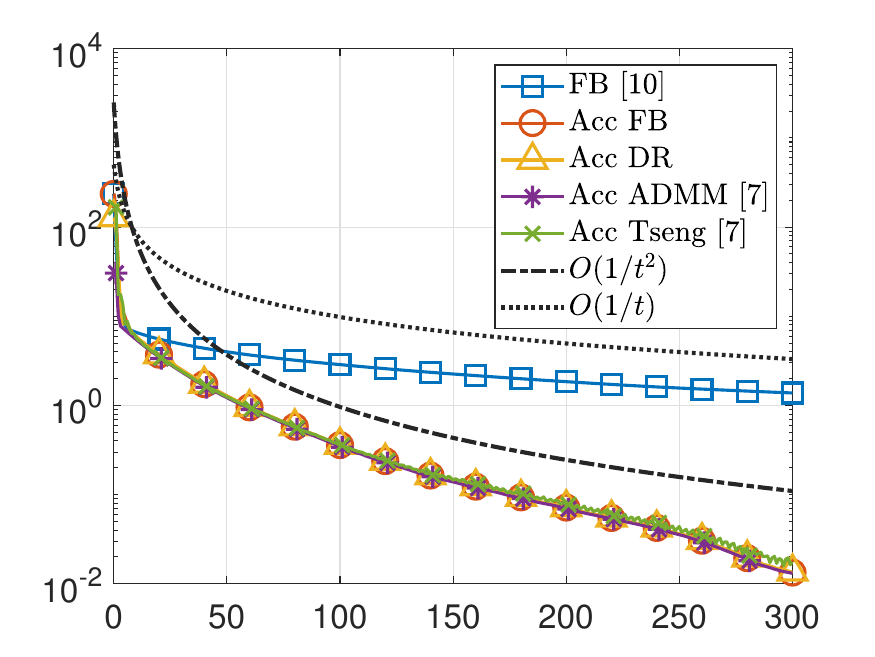}
			\\[-0.3 cm] { $t$}
		\end{tabular}
		&
		\begin{tabular}{c}
			\vspace{.25cm}
			\normalsize{\rotatebox{90}{$\norm{x(t) \,-\, \xs}^2_2/\norm{\xs}^2_2$}}
		\end{tabular}
		&
		\begin{tabular}{c}
			\includegraphics[width=0.28\textwidth]{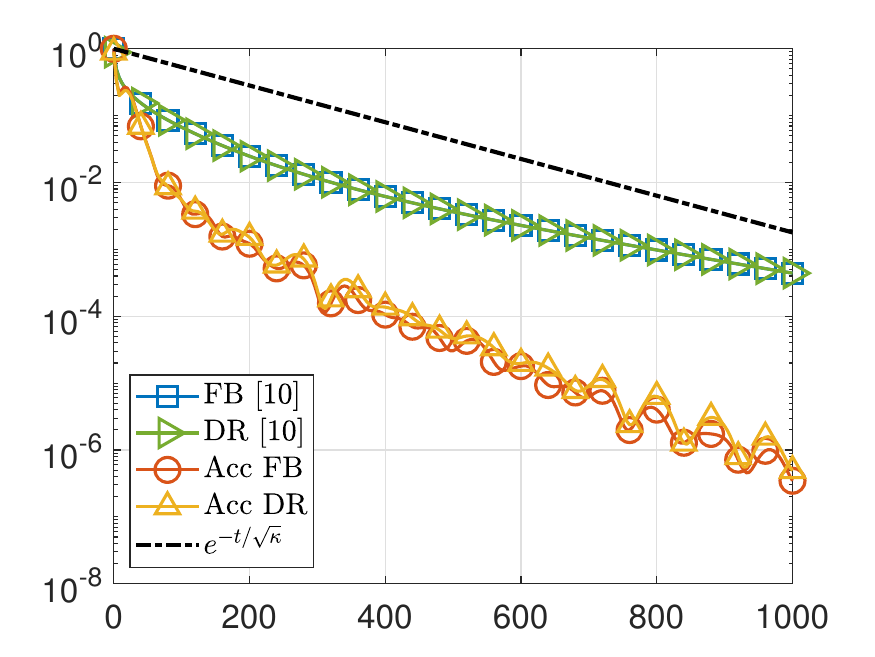}
			\\[-0.3 cm]  {$t$}
		\end{tabular}
		&
		\begin{tabular}{c}
			\vspace{.25cm}
			\normalsize{\rotatebox{90}{$\norm{x(t) \,-\, \xs}^2_2/\norm{\xs}^2_2$}}
		\end{tabular}
		&
		\begin{tabular}{c}
			\includegraphics[width=0.28\textwidth]{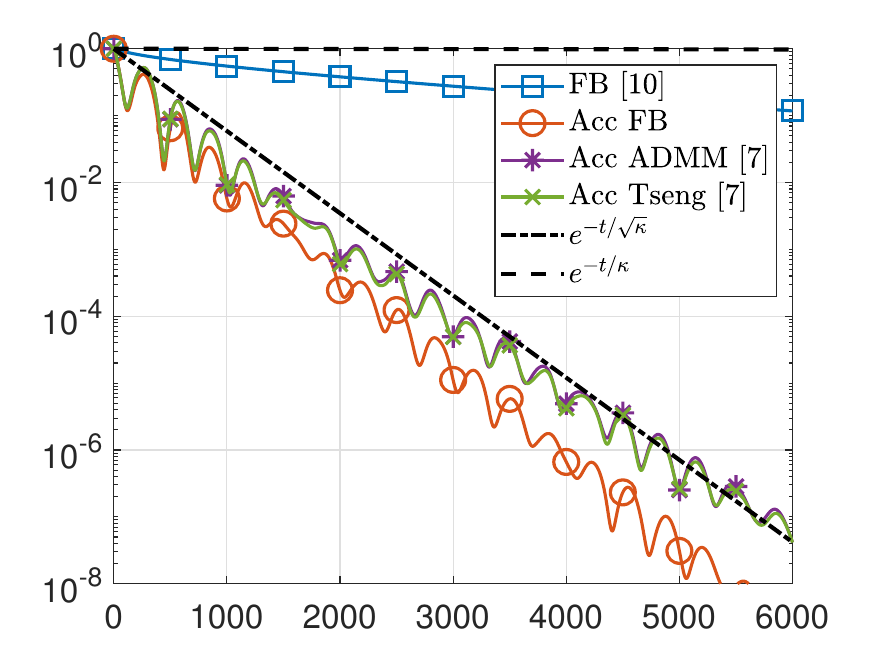}
			\\[-0.3 cm]  {$t$}
		\end{tabular}	
	\end{tabular}
	\vspace{0cm}
	\caption{(a) Sublinear convergence of accelerated dynamics~\eqref{eq.acc_fb} and~\eqref{eq.acc_dr} with the time-varying parameters given in Theorem~\ref{theorem.quad_cvx} for $\ell_1$-regularized least squares problem~\eqref{eq.ls-ell1}. (b) Linear convergence of accelerated dynamics~\eqref{eq.acc_fb} and~\eqref{eq.acc_dr} with the constant parameters given in Theorem~\ref{theorem.quad_str} for box-constrained quadratic program~\eqref{eq.qp-box}. (c) Linear convergence of accelerated dynamics~\eqref{eq.acc_fb} eith the constant parameters given in Theorem~\ref{theorem.gen_str} and $\mu = 1/(L\sqrt[4]{\kappa}) \approx 1/(26L)$ for the $\ell_1$-regularized logistic regression problem~\eqref{eq.l1-logistic} with $\varrho = 0.1$, and condition number $\kappa\approx5\cdot10^5$.}
\end{figure*}

	\vspace*{-1ex}
We examine the $\ell_1$-regularized logistic regression problem
\beq\label{eq.l1-logistic}
\minimize\limits_{x}\sum\limits_{i \,=\, 1}^n\!\!\left( -y_ia_i^Tx  +  \log(1+\rme^{a_i^Tx})\right)    + \lambda\norm{x}_1  +  \tfrac{\varrho}{2}\norm{x}_2^2
\eeq
where the problem data $A = [a_1 \; \cdots ~a_n]^T\in\reals{200\times1000}$ as well as the scalars $\lambda$ and $y_i$ are generated as described in~\cite{suboycan14}. In Fig.~\ref{fig.3}, we validate linear convergence of accelerated proximal gradient flow dynamics~\eqref{eq.acc_fb} with constant parameters given in Theorem~\ref{theorem.gen_str} and $\mu = 1/(L\sqrt[4]{\kappa}) \approx 1/(26L)$ for this strongly convex problem with $\varrho = 0.1$, and the condition number $\kappa\approx5\cdot10^5$.

	\vspace*{-1ex}
\section{Concluding remarks}\label{sec.conclusion}
	\vspace*{-1ex}
	
	We have introduced continuous-time accelerated FB and DR splitting dynamics for composite optimization problems where the objective function is given by the sum of smooth and nonsmooth terms. Motivated by the variable-metric gradient method interpretation of FB and DR splittings, we have leveraged smoothness of associated envelopes and utilized a Lyapunov-based approach to establish convergence rates for both convex and strongly convex problems. 

\appendix 
	\vspace*{-1ex}
\section{Proof of Lemma~\ref{lemma.ineq}}\label{proof.lemma.ineq}
	\vspace*{-1ex}
	
Our proof is based on an equivalent definition of the FB envelope~\cite[Eq. (2.8b)]{patstebem14} 
\beq\label{eq.fb_eq}
\!\cF_{\mu}(x)  =  f(x)  +  g(p(x))  -  \mu\inner{\nabla f(x)}{G_\mu(x)} + \tfrac{\mu}{2}\norm{G_\mu(x)}_2^2
\eeq
and the subgradient inequality,
\beq\label{eq.subgrad}
g(x)\geq g(\xp) + \inner{r}{x-\xp},~\forall x,\xp\in\reals{n},~\forall r\in\partial g(\xp).
\eeq 
\textbf{Proof of~\eqref{eq.lemma.1}.} Let $p(x)\DefinedAs \prox_{\mu g}(x-\mu\nabla f(x))$. The identity~\eqref{eq.moreau_grad} together with~\eqref{eq.gen_grad} implies that $$G_\mu(x) - \nabla f(x) \in\partial g(p(x)).$$ Hence, the subgradient inequality~\eqref{eq.subgrad} for an arbitrary $\xp\in\reals{n}$ yields 
\beq\non
g(\xp) \;\geq\; g(p(x))  \,+\, \inner{G(x) \,-\, \nabla f(x)}{\xp  \,-\,  p(x)}
\eeq
which in conjunction with~\eqref{eq.fb_eq} results in
\beq\non
\ba{rcl}
\cF_{\mu}(x)  &\asp{\leq}&  f(x)   \,+\,   g(\xp)   \,-\,  \inner{G(x)  \,-\,  \nabla f(x)}{\xp  \,-\,   p(x)} 
\\ 
&\asp{}& \,-\,   \mu\inner{\nabla f(x)}{G_\mu(x)}  \,+\,  \frac{\mu}{2}\norm{G_\mu(x)}_2^2.
\ea
\eeq
Rearranging terms and then upper bounding $f(x)$ via strong convexity,
\beq\non
f(x)  \,\leq\, f(\xp)  \,+\, \inner{\nabla f(x)}{x \,-\, \xp} \,-\, \tfrac{m}{2}\norm{x \,-\, \xp}^2
\eeq
 gives
\beq\non
\ba{rcl}
\cF_{\mu}(x)  &\asp{\leq}&  f(\xp)   \,+\,   g(\xp)    \,-\, \frac{m}{2}\norm{x \,-\, \xp}_2^2  
\\ 
&\asp{}& +\,  \inner{G(x)}{p(x)  \,-\,  \xp}  \,+\,   \frac{\mu}{2}\norm{G_\mu(x)}_2^2
\\
&\asp{=}&  F(\xp)    \,-\, \frac{m}{2}\norm{x \,-\, \xp}_2^2 \,+\,  \inner{G_\mu(x)}{x  \,-\,  \xp} 
\\ 
&\asp{}&    \,-\, \inner{G_\mu(x)}{x \,-\, p(x)}    \,+\,  \frac{\mu}{2}\norm{G_\mu(x)}_2^2
\ea
\eeq
where the equality is obtained by adding and subtracting $x$ from the second argument of the inner product. Substituting~\eqref{eq.gen_grad} into the second inner product completes the proof. 

\textbf{Proof of~\eqref{eq.lemma.2}.} Let $p(x)\DefinedAs \prox_{\mu g}(x-\mu\nabla f(x))$, $\xt \DefinedAs x - \xs$, and $\pt \DefinedAs p(x) - p(\xs)$. Substitution of~\eqref{eq.gen_grad} into~\eqref{eq.fb_eq}  gives
\beq\label{eq.fb_eq2}
\ba{rcl}
\cF_\mu(x) \,-\, F(\xs) &\asp{=}& f(x) \,-\, f(\xs) \,+\, g(p(x)) \,-\, g(\xs)
\\ 
&\asp{}& \,-\, \inner{\nabla f(x)}{\xt \,-\, \pt} \,+\, \frac{\mu}{2}\norm{G_\mu(x)}^2.
\ea
\eeq
The optimality condition for problem~\eqref{eq.composite} implies that $-\nabla f(\xs) \in\partial g(\xs)$ and $\xs = p(\xs)$. Using these two identities in the subgradient inequality~\eqref{eq.subgrad} yields
\beq\non
g(p(x)) \,\geq\, g(\xs)  \,-\, \inner{\nabla f(\xs)}{\pt}
\eeq
which together with~\eqref{eq.fb_eq2} results in
\beq\non
\ba{rcl}
\cF_\mu(x) \,-\, F(\xs) &\asp{\geq}& f(x)  \,-\,  f(\xs) \,-\, \inner{\nabla f(x)}{\xt} 
\\ 
&\asp{}& \!\!\!\!+\inner{\nabla f(x)   \,-\,  \nabla f(\xs)}{\pt}  \,+\,  \frac{\mu}{2}\norm{G_\mu(x)}^2.
\ea
\eeq
Substituion of the lower bound given by the cocoercivity of gradient $\nabla f(x)$,
\beq\non
f(x) \,-\, f(\xs)  \;\geq\; \inner{\nabla f(\xs)}{\xt} \,+\, \tfrac{1}{2L}\norm{\nabla f(x)  \,-\, \nabla f(\xs)}_2^2
\eeq
into the right-hand-side yields
\beq\non
\ba{rcl}
\cF_\mu(x)  \,-\,  F(\xs) &\asp{\geq}& \inner{\nabla f(\xs)  \,-\,  \nabla f(x)}{\xt \,-\, \pt} 
\\
&\asp{}&\!\!\!\!+ \tfrac{1}{2L}\norm{\nabla f(x)   \,-\,  \nabla f(\xs)}_2^2 \,+\,  \frac{\mu}{2}\norm{G_\mu(x)}^2
\\
&\asp{=}&\tfrac{\mu}{2}\norm{G_{\mu}(x) - (\nabla f(x) \,-\, \nabla f(\xs))}^2
\\
&&+\, \tfrac{1}{2L}(1 \,-\, \mu L)\norm{\nabla f(x) \,-\, \nabla f(\xs)}^2
\\
&\asp{\geq}& \tfrac{m^2}{2L}(1 \,-\, \mu L)\norm{x \,-\, \xs}_2^2
\ea
\eeq
where the equality follows from~\eqref{eq.gen_grad} and the completion of squares; the second inequality is obtained by ignoring the first (nonnegative) term and using $m$-strong convexity of~$f$ on the second one.

	\vspace*{-2ex}

\end{document}